\documentclass[11pt]{article}
\RequirePackage[OT1]{fontenc}
\RequirePackage{amsthm,amsmath}
\RequirePackage{natbib}
\RequirePackage[colorlinks,citecolor=blue,urlcolor=blue]{hyperref}
\usepackage{epsfig,multicol,multirow}
\usepackage{graphicx}
\usepackage{float}
\usepackage{placeins}
\usepackage{booktabs}
\usepackage{lscape}
\usepackage[table]{xcolor}
\usepackage{adjustbox}
\usepackage{listings}
\usepackage{fullpage}
\usepackage{bm} 

\def\cvf{\mathcal{CV}}
\def\Tf{\mathcal{T}}
\def\Mf{\mathcal{M}}
\def\Vf{\mathcal{V}}
\def\madf{\mathcal{MAD}}
\def\rcvf{\mathcal{RCV}}
\def\Gf{\mathcal{G}}
\def\MADf{\mathcal{MAD}}
\def\e{\hbox{E}}
\def\IF{\hbox{IF}}

\def\var{\hbox{Var}}

\def\iqr{\hbox{IQR}}

\def\cv{\hbox{CV}}
\def\wcv{\widehat{\hbox{CV}}}
\def\scv{\widehat{\hbox{cv}}}
\def\tcv{\widetilde{\hbox{cv}}}
\def\srcv{\widehat{\hbox{rcv}}}
\def\mad{\hbox{MAD}}
\def\rcv{\hbox{RCV}}

\def\ase{\hbox{ASE}}
\def\asd{\hbox{ASD}}
\def\rasd{\hbox{rASD}}
\def\rase{\hbox{rASE}}
\def\asv{\hbox{ASV}}
\def\rase{\hbox{RASE}}

\newcommand{\F}{\mathcal{F}\,}
\newcommand{\asim}{\stackrel{\text{\tiny approx.}}{\sim}}

\numberwithin{equation}{section}
\theoremstyle{plain}
\newtheorem{theorem}{Theorem}[section]

\begin{document}

\title{Robust analogs to the Coefficient of Variation}

\author{Chandima N. P. G. Arachchige \\ Department of Mathematics and Statistics, La Trobe University \\ \href{mailto:18201070@students.latrobe.edu.au}{18201070@students.latrobe.edu.au}\\ \\ Luke A. Prendergast\\ Department of Mathematics and Statistics, La Trobe University \\ \href{mailto:luke.prendergast@latrobe.edu.au}{luke.prendergast@latrobe.edu.au}\\ \\ Robert G. Staudte\\ Department of Mathematics and Statistics, La Trobe University \\ \href{r.staudte@latrobe.edu.au}{r.staudte@latrobe.edu.au}}

\maketitle

\begin{abstract}
The coefficient of variation (CV) is commonly used to measure relative dispersion. However, since it is based on the sample mean and standard deviation, outliers can adversely affect the CV.  Additionally, for skewed distributions the mean and standard deviation do not have natural interpretations and, consequently, neither does the $\cv$.  Here we investigate the extent to which quantile-based measures of relative dispersion can provide appropriate summary information as an alternative to the CV. In particular, we investigate two measures, the first being the interquartile range (in lieu of the standard deviation), divided by the median (in lieu of the mean), and the second being the median absolute deviation (MAD), divided by the median, as robust estimators of relative dispersion. In addition to comparing the influence functions of the competing estimators and their asymptotic biases and variances, we compare interval estimators using simulation studies to assess coverage.
\end{abstract}


\section{Introduction}

The coefficient of variation (CV), defined to be the ratio of the standard deviation to the mean, is the most commonly used method of measuring relative dispersion.  It has applications in many areas, including engineering, physics, chemistry, medicine, economics and finance, to name just a few.   For example, in analytical chemistry the CV is widely used to express the precision and repeatability of an assay \citep{reed2002use}. In finance the coefficient of variation is often considered useful in measuring relative risk \citep{miller1977testing} where a test of the equality of the CVs for two stocks can be performed to compare risk. In economics, the $\cv$ is a summary statistic of inequality \citep[e.g.][]{atkinson1970measurement, chen1996regional}. Other examples use the CV to assess the homogeneity of bone test samples \citep{hamer1995new}, assessing strength of ceramics \citep{gong1999relationship} and as a summary statistic to describe the development of age- and sex-specific cut off points for body-mass indexing in overweight children \citep{cole2000establishing}.   

The lack of robustness to outliers of moment-based measures such as the mean and standard deviation has long been known.  Almost a century ago \cite{lovitt1929statistics} proposed a measure called the \lq\lq coefficient of variability \rq\rq based on the upper and lower quartiles ($Q_3$ and $Q_1$).  Promoted as an alternative to the CV, it was defined to be $(Q_3-Q_1)/(Q_3+Q_1)$.  \cite{bonett2006confidence} have since called this measure the \lq\lq coefficient of quartile variation \rq\rq\  and introduced an interval estimator which exhibited good coverage even for small samples. This measure was recently re-investigated by  \cite{bulent2018bootstrap} and they have constructed bootstrap confidence intervals that typically provide conservative coverage.  Another alternative measure is to take the ratio of the mean absolute deviation from the median divided by the median. This measure has applications in tax assessments \citep{gastwirth1982statistical} and confidence intervals have been considered by  \cite{bonett2005confidence}.  The mean absolute deviation is still non-robust to outliers, and robustness can be improved \cite[see e.g.][]{shapiro2005practical, reimann2008statistical, varmuza2009introduction} by instead using the interquartile range (IQR) or the \textit{median} absolute deviation (MAD).  

For decades, interval estimation for the CV has attracted the attention of many researchers.   For example, \cite{gulhar2012comparison} compared no less than 15 parametric and non-parameic confidence interval estimators of the population CV.  To the best of our knowledge interval estimators have not been introduced for the coefficient of variation based on the IQR and MAD.  Therefore, given the obvious need for interval estimators that has attracted the interest for many others, one aim of this paper is to provide reliable interval estimators.  We are motivated to do so by noting the excellent coverage achieved for measures based on ratios of quantiles, even for small samples \citep[][]{prst-2016b, prst-2017a, prst-2017b, arachchige2019robust}.

\section{Notations and some selected methods}
Let $X_1, X_2, ......, X_n$ be an independent and identically distributed sample of size $n$ from a distribution with distribution function $F$. Then the sample mean estimator is $\overline{X} = n^{-1}\sum^n_iX_i$ and sample variance estimator is $S^2=\sum _i^n(X_i-\overline X)^2/(n-1)$.  The sample coefficient of variation estimator is then $\wcv =S/\overline X$. Next let $\F $ be the class of all right-continuous cdfs on the positive axis; that is each  $F\in \F$ satisfies $F(0)=0.$  For a sample denoted $x_1, \ldots, x_n$, the statistics $\overline x$, $s$, and $\scv=s/\overline{x}$ are the observed values of the $\overline X$, $S$ and $\wcv$ estimators above, and are therefore estimates of the unknown population parameters $\mu =\e _F[X]$, $\sigma =\sqrt {\e _F[(X-\mu )^2]}$ and $\cv =\sigma /\mu $, assuming the first two moments of $F$ exist.

 For each such $F\in \F$ define the associated left-continuous {\em quantile function} of $F$ by $ Q(u)\equiv \inf \{x:\ F(x)\ge u \}$, for $0< u< 1.$  When the population $F$ is understood to be fixed but unknown, we sometimes simply write $x_u=Q(u)$ and write the corresponding estimators of these population quantiles as $\widehat x_{u}$.  We restrict attention to the quartiles $x_{0.25}$, $x_{0.5}$ and $x _{0.75}$, the sample estimates of which we denote $q_1$, $m$ and $q_3$ for convenience.
\subsection{Selected interval estimators of the CV} \label{sec:selected_methods}

We begin by describing the inverse method \citep{sharma1994asymptotic} for obtaining an interval estimator for the $\cv$ since it is perhaps the most naturally arising interval involving only basic principles.  As additional methods for comparison later, we have chosen four of the 15 considered in  \cite{gulhar2012comparison} that exhibited comparatively good performance in terms of coverage. 

While parametric interval estimators for the CV have typically been developed assuming an underlying normal distribution, such as those that we present below, for large sample sizes, they can also perform well \citep{gulhar2012comparison} when there are deviations from normality due to the Central Limit Theorem.

\subsubsection*{The inverse method}\label{sec:Inverse}

Using the above notation, for suitably large $n$, $\overline x/s$ is approximately $N(0,1/n)$ distributed.  An approximate $(1-\alpha/2)\times 100$\% confidence interval for $\mu/\sigma$ is therefore  $\overline x/s \pm z_{1-\alpha/2}/\sqrt{n}$.  Noting that $\mu/\sigma$ is simply the inverse of the population $\cv$, an approximate 95\% confidence interval for the $\cv$ can therefore be obtained by inverting this interval for $\mu/\sigma$, giving \citep{sharma1994asymptotic}
\begin{align}\label{eq_inverseMethod}
\left\{\left[\frac{1}{\scv} \ + \ z_{1-\alpha/2}\left(\frac{1}{n^{1/2}}\right)\right]^{-1} ,\  \left[\frac{1}{\scv} \ - \ z_{1-\alpha/2} \left(\frac{1}{n^{1/2}}\right)\right]^{-1} \right\}~.
\end{align}
Robustness of this interval estimator was recently re-investigated by \cite{groen-2011}.

\subsubsection*{The median-modified Miller interval (Med Mill)}\label{sec:Med Mill}
The $\cv$ estimator has an approximate asymptotic normal distribution with mean $\cv$ and variance $(n-1)^{-1}\cv^2(0.5+\cv^2)$ leading to an asymptotic interval proposed by  \cite{edward1991asymptotic}. In noting that the mean is a poor summary statistic of central location for skewed distributions,  \cite{gulhar2012comparison} proposed a median modification where the sample median replaces the sample mean in $s$.  Let $\tilde{s} = \sqrt{\frac{1}{n-1} \sum_{i=1}^{n}(x_i-  m)^2}$ and 
$\tcv=\tilde{s}/\overline{x}$, the interval estimator is 
\begin{equation} \label{eq:MedMill}
\resizebox{\textwidth}{!}{$\left\{\tcv - z_{1-\alpha/2}\sqrt{(n-1)^{-1}\tcv^2\left(0.5+\tcv^2\right)}, \ \tcv + z_{1-\alpha/2}\sqrt{(n-1)^{-1}\tcv^2\left(0.5+\tcv^2\right)}\right\}$}~.
\end{equation}
While simulations conducted by \cite{gulhar2012comparison} using data sampled from a chi-square and gamma distribution showed typically good results for the  \cite{edward1991asymptotic} interval, coverage was often better, if not at least similar, when using the median modification.  With our interest mainly in skewed distributions, we focus on the median modified interval in \eqref{eq:MedMill}. 

\subsubsection*{Median modification of the modified McKay (Med MMcK)}\label{sec:Med MMcK}
 \cite{gulhar2012comparison} also introduced a median modification to the \textit{modified McKay interval} \citep{mckay1932distribution, vangel1996confidence}.  The median-modified interval is
\begin{equation}\label{eq:MedMMcK}
\resizebox{\textwidth}{!}{$\left\{\tcv\sqrt{\left(\frac{ \chi_{n-1, 1-\alpha/2}^2 +2}{n}-1\right)\tcv^2 + \frac{ \chi_{n-1,1-\alpha/2}^2}{n-1}},\ 
\tcv\sqrt{\left(\frac{ \chi_{n-1,\alpha/2}^2 +2}{n}-1\right)\tcv^2 + \frac{ \chi_{n-1,\alpha/2}^2}{n-1}}\right\}$} ~,
\end{equation}
where $\chi_{n-1,\alpha}^2$  is  the  $100\alpha$-th  percentile  of a  chi-square   distribution with $(n-1)$ degrees of freedom. We focus on this median modified interval based on the results in \cite{gulhar2012comparison}. 

\subsubsection*{The Panich method}\label{sec:Panich}
 \cite{panichkitkosolkul2009improved} has further modified the Modified McKay \citep{vangel1996confidence} interval by replacing the sample $\cv$ with the maximum likelihood estimator for a normal distribution, $\tilde{k}=\sqrt{\sum_{i=1}^{n}(x_i- \overline{x})^2}/(\sqrt{n}\overline{x})$.  The interval is 
\begin{equation}\label{eq:Panich}
\resizebox{\textwidth}{!}{$\left\{\tilde{k}\sqrt{\left(\frac{ \chi_{n-1, 1-\alpha/2}^2 +2}{n}-1\right)\tilde{k}^2 + \frac{ \chi_{n-1,1-\alpha/2}^2}{n-1}},\ 
\tilde{k}\sqrt{\left(\frac{ \chi_{n-1,\alpha/2}^2 +2}{n}-1\right)\tilde{k}^2 + \frac{ \chi_{n-1,\alpha/2}^2}{n-1}}\right\}$}~.
\end{equation}

\subsubsection*{The Gulhar method}\label{sec:Gulhar}
Using the fact that $(n-1)S^2/\sigma^2\sim \chi^2_{n-1}$ when data is sampled from the normal distribution,  \cite{gulhar2012comparison} proposed the interval,
\begin{equation}\label{eq:Gulhar}
\left(\frac{\sqrt{(n-1)}\,\scv}{\sqrt{\chi_{n-1,1-\alpha/2}^2}}, \ \frac{\sqrt{(n-1)}\,\scv}{\sqrt{\chi_{n-1,\alpha/2}^2}}\right),
\end{equation}
which compared favorably to the median-modified intervals for larger CV values.  We therefore use this interval as one of the competitors.

\subsection{Two robust versions of the CV}
We now consider two robust alternatives for the CV that are based on quantiles. The denominator for the measures is the median, a preferred measure of centrality than the mean for skewed distributions. 

\subsubsection{A version based on the IQR}\label{sec:RCV_Q}
An option for the numerator is to use the interquartile range ($\iqr$). \cite{shapiro2005practical} gives this alternative as
\begin{equation} \label{eq:RCV_Q}
\rcv_Q=0.75 \times \frac{\iqr}{m}~,   
\end{equation}
where the multiplicative factor 0.75 makes $\rcv_Q$ comparable to the CV for a normal distribution.  To the best of our knowledge there has been no research into interval estimators of the $\rcv_Q$ and this will be one of our foci shortly.
 
\subsubsection{A version based on the median absolute deviation}\label{sec:RCV_M}
The median absolute deviation \citep[MAD]{hamp-1974} is defined to be
\begin{equation}\label{eq:MAD}
 \mad =\text{med}\mid{x_i-m}\mid~,  
\end{equation}
where, for \lq $\text{med}$\rq denoting median and i=1,\ldots,$n$. Using the MAD for relative dispersion has been recently proposed \citep[e.g.][]{reimann2008statistical, varmuza2009introduction} giving
\begin{equation}\label{eq:RCV_M}
\rcv_M =1.4826 \times \frac{\mad}{m}~.  
\end{equation}
The multiplier $1.4826=1/\Phi^{-1}(3/4)$, where $\Phi^{-1}$ denotes the quantile function for the $N(0,1)$ distribution, is used to achieve equivalence between $1.4826 \times \mad/m$ and the standard deviation at the normal model. $1.4826 \times \mad/m$ is commonly called the \textit{standardized MAD}. 

\section{Some comparisons between the measures}

 The question of interest is, can we do just as well (or better) in assessing the relative dispersion by replacing the population concepts
 $\mu $ and $\sigma $  by the median $m=x_{0.5}$ and interquartile range $\iqr = q_3 - q_1$ or the MAD?

\begin{table}[h!t]
  \centering
  \scriptsize
  \caption{A comparison of the $\cv$, $\rcv_Q$ and $\rcv_M$ for several distributions. LN refers to the log-normal distribution, WEI$(\lambda,\alpha)$ and PAR$(\lambda,\alpha)$ to the Weibull and Pareto Type II distributions with scale parameter $\lambda$ and shape parameter $\alpha$.}
  \resizebox{\linewidth}{!}{%
    \begin{tabular}{cccc}
    \toprule
    Distribution & \multicolumn{1}{c}{$\cv$} & \multicolumn{1}{c}{0.75* IQR/$m$} & \multicolumn{1}{c}{1.4826*MAD/$m$}\\
    \midrule
    Normal($\mu$,$\sigma^2$) &  $\displaystyle\frac{\sigma}{\mu}$    &  $\displaystyle\frac{3}{4}\frac{\sigma}{\mu}\left[\Phi^{-1}(0.75)-\Phi^{-1}(0.25)\right]$      &  $\displaystyle\frac{\sigma}{\mu}$\\
    EXP($\lambda$)   & 1     &   1.189    & 1.030 \\
  Uniform$(a,b)$  & $\displaystyle\frac{1}{\sqrt{3}}\cdot\frac{(b - a)}{(b + a)}$ & $\displaystyle\frac{3}{4}\cdot\frac{(b - a)}{(b + a)}$ & $\displaystyle\frac{1}{\Phi^{-1}(3/4)}\cdot\frac{(b - a)}{(b + a)}$ \\
    WEI($\lambda$, 1) & 1 & 1.189 & 1.029\\
    WEI($\lambda$, 2) & 0.523 & 0.578 & 0.565\\
    WEI($\lambda$, 5) & 0.229 & 0.232 & 0.229\\
    $\chi^2_2$ & 1 & 1.189 & 1.030 \\
    $\chi^2_5$ & 0.632 & 0.681 & 0.646 \\
    $\chi^2_{\nu \rightarrow \infty}$ & $\rightarrow 0$ & $\rightarrow 0$ & $\rightarrow 0$\\
    LN$(\mu, 1)$ & 1.311 & 1.090 & 0.888 \\
    LN$(\mu, 2)$ & 7.321 & 2.695 & 1.333  \\
    PAR$(\lambda, 2.5)$ & 2.236 & 1.453 & 1.120 \\
    PAR$(\lambda, 5)$ & 1.291 & 1.313 &  1.077\\
    \bottomrule
    \end{tabular}}%
  \label{table:compare}%
\end{table}%

In Table \ref{table:compare} we compare the $\cv$, $\rcv_Q$ and $\rcv_M$ for several distributions. In most cases, the results show an approximate equivalence between the three measures when the underlying population is normal and closer agreement between the two for many other distributions. Hereafter our main interest is comparing the concepts
$\cv$, $\rcv_Q$ and $\rcv_M$ and the natural estimators of them.  

\subsection{Properties}\label{sec:Properties}

An essential property of a measure of relative dispersion is scale invariance. The $\cv $ is well-established, so competing measures should give roughly the same values when the underlying distribution is uni-modal and skewed to the right, As we have seen by examples, the plug-in estimator $s/\bar x$ of $\cv $ suffers from over-sensitivity to outliers. Table \ref{table:prope} provides a rough summary of results in this work.

\begin{table}[ht]
\begin{center}
\scriptsize
\caption{Desirable properties of measures of dispersion and their estimators.  Here \lq +\rq,\  \lq 0\rq\  and \lq $-$\rq\  indicate the property always, sometimes or never holds.}
\vspace{0.5cm}
 \resizebox{\linewidth}{!}{%
\begin{tabular}{lccccc}
   Property  &&   $\cv$   & $\rcv _Q$ & $\rcv _M$ \\
  \hline
 P1: Scale invariant             && +    & + & +             \\
 P2: Simple to understand        && +    & + & 0            \\
 P3: Widely accepted and used    && +    & $0$ & 0       \\
 P4: Defined for all $F$         && $0$\footnotemark[1]  & + & +           \\
 P5: Bounded influence function  && $-$  & + & +          \\[.5cm]
   Property  &&   $\widehat{\cv} $   & $\widehat{\rcv} _Q$ & $\widehat{\rcv} _M$\\
 \hline
 P6: Consistency                                 &&  0\footnotemark[2] & + & +   \\
 P7: Asymptotic normality                        &&  0   & +  & +  \\
 P8: Standard error formula available            &&  +   & + & +    \\
 P9: Unaffected by 1\% moderate outliers         &&  0   & + & +  \\
 P10: Unaffected by 1\% extreme outliers         &&  $-$ & + & +    \\
 P11: Reliable coverage of confidence intervals  &&  $-$ & + & \ +
 \label{table:prope}
\end{tabular}}
\end{center}
\end{table}
\FloatBarrier
\vspace{-1cm}
\footnotetext[1]{The $\cv $ is only defined if $F$ has a finite variance, but this is usually satisfied for
diameter distribution models.}
\footnotetext[2]{Consistency and asymptotic normality require the existence of certain moments for $F$.}
\bigskip
In the next section, we briefly describe the methodology required to find standard errors and confidence intervals for $\cv$, $\rcv_Q$ and $\rcv_M$.
We also investigate the robustness properties of the point estimators using theoretical methods and simulation studies and we illustrate our methods on a real data set. Finally, a summary and discussion of further possible work is in Section~\ref{sec:summary}.

\subsection{Influence functions}\label{sect:influence_functions}

Consider a distribution function $F$ and suppose that a parameter of interest from $F$ is $\theta$.  Let $\Tf$ be a statistical function for estimator of $\theta$ such that $\Tf(F)=\theta$ and $\Tf(F_n)=\widehat{\theta}$, for $F_n$ denoting an empirical distribution function for sample of $n$ observations from $F$, denotes an estimate of $\theta$.  Now, for $0\leq \epsilon \leq 1$, define the \lq contamination\rq distribution $(F_\epsilon)$ to have positive  probability $\epsilon $ on $x$ (the contamination point) and $1-\epsilon $ on the distribution $F$ such that $F_\epsilon=(1-\epsilon)F+\epsilon\Delta_x$ where $\Delta_x$ denotes the distribution function that puts all of its mass at the point $x$. The influence of the contamination on the estimator with functional $\Tf$, relative to proportion of contamination, is $[\Tf(F_\epsilon) - \Tf(F)]/\epsilon$.
The influence function  \citep{hamp-1974} is then defined for each $x$ as
  $$\IF (x;\Tf,F)=\lim_{\epsilon\downarrow 0} \frac{\Tf(F_\epsilon) - \Tf(F)}{\epsilon}\equiv \frac {\partial }{\partial \epsilon }\Tf(F_\epsilon)\Big|_{\epsilon = 0}.$$

A convenient way to appreciate the usefulness of the influence function in studying estimators is to consider the power series expansion $\Tf(F_\epsilon)=T(F)+\epsilon \IF (x;\Tf,F) + O(\epsilon^2)$.  So that, ignoring the error term $O(\epsilon^2)$ which is negligible for small $\epsilon$, increasing $\left|\IF (x;\Tf,F)\right|$ results in increasing influence of contamination on the estimator. Consequently, the influence function provides a very useful tool in the study of robustness of estimators.

One can show that \citep[e.g.,][]{HRRS86,S-S-1990} for $X\sim F$, the mean and variance at $F$ of the random influence function are $\e _F [\IF (X;\,\Tf,F)]=0$ and $\var _F [\IF (X;\,\Tf,F)]=\e _F [\IF ^2(X;\,\Tf,F)]$. A reason for finding this last variance is that it arises in the asymptotic variance of the functional  of $\Tf(F_n)$; that is, 
\begin{equation}
n\;\var [\Tf(F_n)]\to \asv\left(\Tf,F\right) = \e _F [\IF ^2(X;\,\Tf,F)]~.\label{asv}
\end{equation}

\subsubsection{Influence function of the CV} \label{sec:IF_cv}
Let $\Mf$ and $\Vf$ denote the functional for the usual mean and variance estimators such that, at $F$, $\Mf(F)=\int xdF=\mu$ and $\Vf(F)=\int\left[x-\Mf(F)\right]^2df=\sigma^2$. The respective influence functions are
$\IF (x;\,\Mf,F)=x-\mu$ and $\IF (x;\,\Vf,F)=(x-\mu)^2-\sigma^2$.   For convenience in notation, let $\cvf$ also denote the functional for the $\cv$. \cite{groen-2011} derives the influence function as
\begin{equation}\label{IFcv}
  \IF (x;\cvf ,F)=\cv\left [\frac{\IF (x;\,\Vf,F)}{2\sigma^2}-\frac{\IF (x;\,\Mf,F)}{\mu} \right].
\end{equation}

\subsubsection{Influence function of the IQR-based RCV}\label{sec:IF_RCV_Q}
The influence function of the $p$th quantile $x_p=\Gf(F;p)=F^{-1}(p)$ is well-known \cite[p.59]{S-S-1990} to be
$\IF [x;\,\Gf(\,\cdot ,p),F]= \{p-I[x_p\geq x]\}\,g(p)$,
where $\Gf'(F;p)=g(p)= 1/f(x_p)$ is the {\em quantile density} of $\Gf$ at $p$.
The influence function of the ratio of two quantiles $\rho _{p,q}(F)=x_p/x_q=\Gf(\,\cdot  ,p)/\Gf(\,\cdot ,q)$ is then found to be \cite{prst-2017a}:
\begin{equation}\label{eq:IF_rho}
  \IF (x;\,\rho _{p,q} ,F) = \rho _{p,q} \left\{\frac {\IF [x;\,\Gf(\,\cdot ,p),F]}{x_{p}}
         -\frac {\IF [x;\,\Gf(\,\cdot ,q),F]}{x_{q}}\right\}~.
 \end{equation}
It then follows that the influence function of $\rcvf_Q(F)=0.75\,\iqr /m$ in terms of  (\ref{eq:IF_rho}) is
\begin{equation}\label{eq:IF_RCVq}
  \IF (x;\,\rcvf_Q ,F) = 0.75\,\left[ \IF (x;\,\rho _{3/4,1/2} ,F)-\IF (x;\,\rho _{1/4,1/2} ,F)\right ]~.
 \end{equation}
 
\subsubsection{Influence function of the MAD-based RCV}\label{sec:IF_RCV_M}
Let $\MADf$ denote the functional for the standardized MAD.  The influence function for the MAD estimator was described by \cite{hamp-1974} and its form for the standardized MAD for the standard normal distribution is \citep[see, e.g., page 107 of][]{HRRS86}
\begin{equation}
\IF (x;\,\MADf ,\Phi) = \frac{1}{4\Phi^{-1}(0.75)\phi\left[\Phi^{-1}(0.75)\right]}\text{sign}\left[|x| - \Phi^{-1}(0.75)\right]~.\label{IFmad}   
\end{equation}

It is not suitable for us to study the influence function for $\rcvf_M$ at the standard normal model since the median is equal to zero.  However, the influence function for the standardized MAD for an arbitrary mean, $\mu$, for the normal distribution is simply \eqref{IFmad} shifted to be centred at $\mu$ and therefore equal to
$\IF (x;\,\MADf ,\Phi_\mu)= \IF (x-\mu;\,\MADf ,\Phi)$ where we let $\Phi_\mu$ denote the distribution function for the $N(\mu, 1)$ distribution.

Let $\rcvf_M$ be the statistical functional for the MAD-based RCV such that $\rcvf_M(F)=\MADf(F)/\Gf(F,1/2)=\rcv_M$.  Hence, using the Product Rule and the Chain Rule, the influence function for the RCV$_M$ estimator is
\begin{align}
 \text{IF}(x;\rcvf_M, F)&=\frac{\partial}{\partial\epsilon}\rcvf_M(F_\epsilon^{(x)})\big|_{\epsilon=0} \nonumber\\
 &=\frac{\IF (x;\,\MADf ,\Phi_\mu)}{m}-\text{RCV}_M\frac{\IF (x;\,\Gf(\,\cdot ,1/2),F)}{m}~.   
\end{align}  

The general form of the influence for the MAD can be found in, for example, page 137 of  \cite{huber1981robust}, page 16 of \cite{andersen2008modern} and page 37 of  \cite{wilcox2011introduction} and this will be used to plot the influence functions for the non-Gaussian examples that follow.

\subsubsection{Example influence function comparisons}
 To compute the true value for the MAD for the distributions being considered for influence function comparisons, and also when required later, we used the R function we have provided in Section \ref{Rcode}. Readers can use this code to compute the true MAD for any distributions.  
  
\begin{figure}
    \centering
    \includegraphics[width=\linewidth]{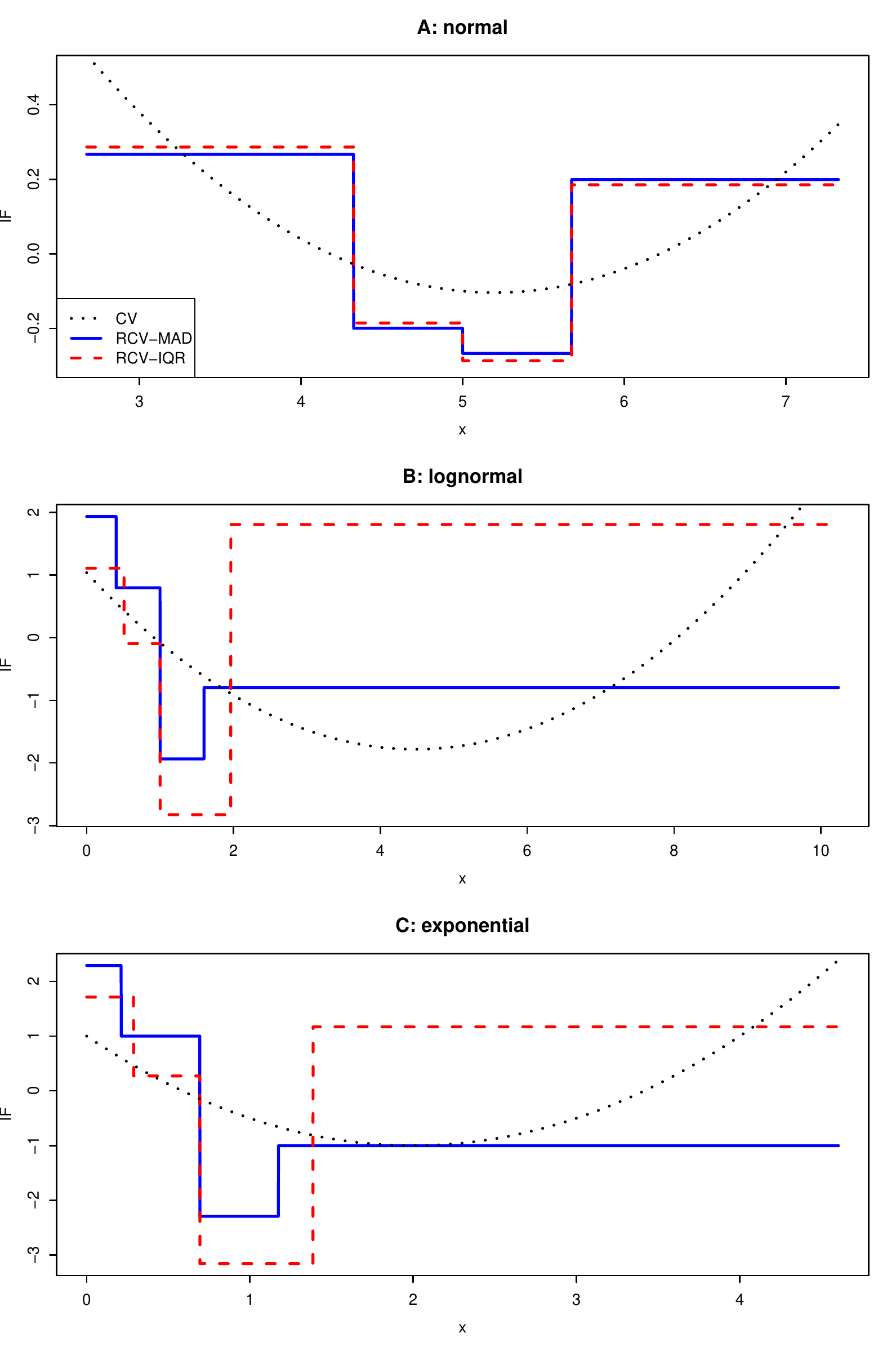}
    \caption{Influence function comparisons between the three measures: CV (black, solid), RCV$_M$ (blue, dash) and RCV$_Q$ (red, dots) for (A) the normal, (B) log-normal and (C) exponential distributions.}
    \label{fig:IFplot}
\end{figure}

In Plot A of Figure \ref{fig:IFplot} we plot the influence functions for the three measures.  The influence functions for the two robust measures are almost identical.  In fact, it is know that the influence functions for the IQR and MAD are the same for the normal distribution \citep[see page 110 of][]{HRRS86} so that the measures share the same robustness properties for this model. The differences in Figure \ref{fig:IFplot} are due to the multiplier 0.75 for the IQR based measure chosen to give approximate equivalence, instead of exact, for the normal. However, this does not generalize to all distributions.  As expected, the influence function for the CV is unbounded, meaning that outliers are expected to have uncapped influence on the estimator as they move further from the population mean.  On the other hand, the influence functions for the robust measures are bounded.  Extreme outliers are expected to have no more influence on the estimators when compared to, say, those closer to the 25\% and 75\% percentiles.  However, the discontinuities at the median and the 25\% and 75\% percentiles, suggest that the estimators are more sensitive locally in these areas. 

\subsection{Asymptotic variances and standard deviations}\label{sec:ASV}
In this section, we further compare the estimators by deriving their asymptotic variances. As discussed in Section \ref{sect:influence_functions}, for an estimator with functional $\Tf$, the asymptotic standard deviation can be found by $\asd (\Tf,F)\equiv \sqrt{\asv(\Tf,F)}=\sqrt{\{\e _F [\IF^2 (X;\,\Tf,F)]\}}$.  We now derive the ASVs for the estimators before comparing their relative asymptotic standard deviations. 

\subsubsection{Asymptotic Variance of the CV estimator}\label{sec:ASV_CV}

Recall $\mu =\Mf(F)$ is the mean for distribution $F$ and let $\mu _k= \e_F [\{X-\Mf(F)\}^k]$ denotes the $k$th central moment of $X\sim F$ where $\mu_2=\sigma^2=\Vf(F)$ denotes the variance.  The influence function for the mean is $\IF(x;\Mf,F)=x-\mu$ and $E\left[\IF(X;\Mf,F)^2\right]=\sigma^2=\asv(\Mf,F)$, the asymptotic variance of the mean estimator.  Similarly, $\IF(x;\Vf,F)=(x-\mu)^2-\sigma^2$ and $E\left[ \IF(X;\Vf,F)^2\right]=\mu_4-\sigma^4=\asv(\Vf,F)$.  Before deriving the ASV for the CV estimator, we note that $E\left[\IF(X;\Mf,F)\IF(X;\Vf,F)\right]$, which is the asymptotic covariance between the mean and variance estimators, is equal to $\mu_3-\sigma^2$.  Now, from \eqref{IFcv},

\begin{align}
 \e\left[\IF(X;\cvf,F)^2\right]=&\left[\cvf(F)\right]^2\Bigg\{\frac{\asv(\Vf,F)}{4\sigma^4}+\frac{\asv(\Mf,F)}{\mu^2}\nonumber\\
& \qquad \qquad \quad -\frac{E\left[ \IF(X;\Vf,F)\IF(X;\Mf,F)\right]}{\sigma^2\mu}\Bigg\} \nonumber\\
  \asv(\cvf,F) =&\cv^2\left(\frac{\mu _4-\sigma^4}{4\sigma^4}+\frac{\sigma^2}{\mu^2}-\frac {\mu _3}{\sigma^2\mu}\right)~,
  \label{eqn:asymvarCV}
\end{align}    
assuming that the fourth moment exists.

Note that for $X\sim F$, $\mu_3=0$ and $\mu_4=3\sigma^4$ so that $\asv(\cvf,F)=\cv^2\left(1/2 + \cv^2\right)$ which is the asymptotic variance used by  \cite{edward1991asymptotic} in the construction of the asymptotic interval for the CV detailed in Section \ref{sec:selected_methods}.

\subsubsection{Asymptotic Variance of the $\rcv_Q$ estimator}\label{sec:ASV_RCV_Q}
The asymptotic variance of the estimator of $x_p$, the $p$-th quantile, is well known to be \citep[eg. Ch.2 of ][Ch.3]{david-1981,Das-2008} $\asv\left(\Gf, F; p\right)=p(1-p)g^2(p)$ where, as denoted earlier, $g(p)=1/f(x_p)$ and $f$ is the density function.  This can be verified also using $\e\left[\IF(X;\Gf(\cdot,p), F)^2\right]$.  Similarly, and as also found in the preceding references, the asymptotic covariance between the $p$-th and $q$-th quantile estimators is, $\e\left[\IF(X;\Gf(\cdot,p), F)\IF(X;\Gf(\cdot,q), F)\right]=p(1-q)g(p)g(q)$, provided $0< p < q < 1$.

Asymptotic variance for $\rcv_Q=0.75\,\iqr /m$ is obtained by a straightforward but lengthy derivation of $\e\left[\IF(X;\rcvf_Q,F)^2\right]$ with $\IF(X;\rcvf_Q,F)$ defined in \eqref{eq:IF_RCVq} (or by using the Delta method).  After simplifying, it is

\begin{theorem}\label{th:ASV_RCVQ}
The asymptotic variance for the estimator of $\rcv_Q$ is
\begin{align*}
   \asv(\rcvf_Q,F) =&  
   \frac{\rcv_Q^2}{4}\Bigg\{\frac {3\left[g^2(3/4)+g^2(1/4)\right]-2\, g(3/4)g(1/4)}{4\times IQR^2}  \nonumber \\
   & \qquad \qquad +  \frac {g^2(1/2)}{m^2} - \frac {g(1/2)\left[g(3/4)-g(1/4)\right]} {m\times IQR} \Bigg\}~.  
\end{align*}
\end{theorem}
\noindent The proof of Theorem \ref{th:ASV_RCVQ} is in Section \ref{app:proof_of_ASV_RCVQ}.

\subsubsection{Asymptotic Variance of the $\rcv_M$ estimator}\label{sec:ASV_RCV_M}
 \cite{falk1997asymptotic} proves the asymptotic joint normality of the $m(F_n)$ and $\madf(F_n)$ estimators. Let $f=F'$ be the density function associated with $F$. If $F$ is continuous near and differentiable at $F^{-1}(1/2)$, $F^{-1}(1/2)-\mad$ and $F^{-1}(1/2) + \mad$ with $f(F^{-1}(1/2))>0$ and $C1 =f(F^{-1}(1/2) - \mad) + f(F^{-1}(1/2) + \mad) > 0$, then 
$$ \sqrt{n}[m(F_n)-F^{-1}(1/2),\, \madf(F_n)-\madf(F)]^\top \asim  N(\mathbf{0}, \bm{\Sigma})~,$$
where \lq $\asim$\rq denotes \lq approximately distributed as for suitably large $n$\rq, $\mathbf{0}$ is a column vector zeroes and $\bm{\Sigma}$ is a two-dimensional covariance matrix with $\text{vec}(\bm{\Sigma})=[\rho_1,\rho_{12},\rho_{12},\rho_2]$.  Hence,
$\rho_1$, $\rho_2$ are the asymptotic variances of the median and MAD estimators respectively and $\rho_{12}$ is the asymptotic covariance between the two.  They are \citep[e.g.][]{falk1997asymptotic}, 
$$\displaystyle\rho_1 = \frac{1}{4f^2(F^{-1}(1/2))}\, , \, \displaystyle\rho_2=\frac{1}{4C_1^2}\left[1+\frac{C_2}{\left[f(F^{-1}(1/2)\right]^2}\right]$$  $$\text{and}\,\,\, \displaystyle\rho_{12}=\frac{1}{4C_1f(F^{-1}(1/2))}\left[1-4F(F^{-1}(1/2)-\mad)+ \frac{C_3}{f(F^{-1}(1/2))}\right]$$
where $C_3 =f(F^{-1}(1/2)-\mad)-f(F^{-1}(1/2)+\mad)$ and $C_2=C_3^2 +4C_3f(F^{-1}(1/2))(1-F(F^{-1}(1/2)+\mad)-F(F^{-1}(1/2)-\mad))$~.

Using the above results and the Delta method \citep[see e.g.][]{Das-2008}, we derived the asymptotic variance of the $\rcv_M$ as given below,
\begin{equation}\label{eq:ASV_RCVM}
   \asv(\rcvf_M,F) = \rcv_M^2 \left(\frac{\rho_1}{m^2}+\frac{\rho_2}{\mad^2}-  \frac{ 2\rho_{12}}{ m \times \mad}\right)~ .
\end{equation}

\subsubsection{Relative asymptotic standard deviation comparisons}\label{sec:ASE Comp}
As an example, the asymptotic standard deviation ($\asd$) for the $\rcv_M$ estimator is given as $\asd(\rcvf_M,F)=\sqrt{\asv(\rcvf_M,F)}$ and the $\asd$s for the other estimators are determined similarly.  Later, we will construct approximate confidence intervals for the measures and therefore it make sense that we use the $\ase$ for comparisons here. 
Since the $\cv$, $\rcv_Q$ and $\rcv_M$ represent different values we use the relative (to the population parameter) $\asd$ ($\rase$) to compare the estimators.  For example, for the $\rcv_M$ estimator this is defined to be $\rasd(\rcvf_M,F) = \asd(\rcvf_M,F)/\rcvf_M(F)$.

\begin{table}[H]
  \centering
  \scriptsize
  \caption{Relative $\asd$ ($\rasd$) comparisons for the estimators of $\cv$, $\rcv_Q$ and $\rcv_M$ for the N$(5, \sigma^2)$, LN$(0, \sigma)$,  EXP($\lambda$) and PAR$(\alpha)$ distributions.} 
   \resizebox{\linewidth}{!}{%
    \begin{tabular}{llccc}
    \toprule
   Distribution & & $\rasd$ for the  & $\rasd$ for the & $\rasd$ for the  \\
        &  &  $\cv$ estimator & $\rcv_Q$ estimator & $\rcv_M$ estimator\\
        \midrule
    N$(5, \sigma^2)$& $\sigma = 0.50$   &  0.714 & 1.173 & 1.173 \\
     &$\sigma = 1$  & 0.735   & 1.193 & 1.193 \\
     & $\sigma = 1.5$ &  0.768  & 1.225 & 1.225 \\
     &$\sigma = 2$ & 0.812 & 1.270  &  1.270 \\
     &   $\sigma = 2.5$ & 0.866  & 1.324 & 1.324 \\
     &$\sigma = 3$ & 0.927 & 1.388 &  1.388\\
    \midrule
    LN$(0, \sigma)$ & $\sigma = 0.10$   & 0.721  & 1.172  & 1.164 \\
     &$\sigma = 0.25$  & 0.801  &1.199& 1.149\\
     & $\sigma = 0.5$ & 1.151  & 1.294  & 1.098 \\
     &$\sigma = 0.75$ & 2.075  & 1.438 & 1.017 \\
     &   $\sigma = 1$ & 4.674 & 1.621 & 0.914 \\
     &$\sigma = 1.5$ & 49.298 & 2.062 & 0.669 \\
    \midrule
    EXP$(\lambda)$ & $\lambda$ & 1  & 1.594 & 0.950 \\
     \midrule
    PAR($\alpha$) & $\alpha = 0.50$  & Undefined     & 3.223 & 0.419 \\
     &$\alpha = 1$  & Undefined    & 2.236 & 0.664\\
     & $\alpha = 1.5$ & Undefined     & 1.976 & 0.735 \\
     &$\alpha = 2$ & Undefined     & 1.862 & 0.785\\
     &   $\alpha = 2.5$ & Undefined     & 1.799 & 0.816\\
     &$\alpha = 3$ & Undefined     & 1.760 &  0.837\\
   &$\alpha = 4$   & 54.482 & 1.714 & 0.864  \\
    & $\alpha = 4.5$  & 5.619 & 1.699 & 0.873 \\
    &  $\alpha = 5$ & 3.724 & 1.687 & 0.880\\
    & $\alpha= 5.5$  & 2.937  & 1.678 & 0.887 \\
    &    $\alpha = 6$   & 2.500  & 1.670 & 0.892  \\
    & $\alpha = 6.5$  & 2.221  & 1.664 & 0.897\\
     \bottomrule
    \end{tabular}}%
  \label{table:Rel_SE}%
\end{table}%
 
 To compare the  $\rasd$ for the estimators of $\cv$, $\rcv_Q$ and $\rcv_M$, we have selected normal and lognormal distributions, both with varying $\sigma$, exponential and the Pareto type II distribution with varying shape. From Table \ref{table:Rel_SE}, the $\rasd$ for $\rcv_Q$ and $\rcv_M$ are a little higher than the $\rasd$ of $\cv$ for the normal distribution. However, $\rcv_Q$ and $\rcv_M$ estimators compare favorably to the $\cv$ for skewed distributions such as the lognormal and Pareto. The $p^{th}$ central moment of Pareto type II distribution exists only if $\alpha > p$ so that the rASD for the $\cv$ estimator is undefined for $\alpha < 4$ since it requires the fourth central moment. When comparing $\rcv_q$ and $\rcv_M$, the $\rcv_M$ estimator is the better performer with smaller (or equal to in the case of the normal) $\rasd$.
 
\section{Inference}\label{sec:inference}

We want to compare point and interval estimators of $\cv =\sigma /\mu $ , $\rcv_Q =0.75\,\iqr /x_{0.5}$ and $\rcv_M =1.4826\,\mad /x_{0.5}$. First, we introduce asymptotic Wald-type intervals using the asymptotic standard errors from earlier.  With recent results highlighting very good coverage for estimators based on ratios of quantiles even for small samples \citep[][]{prst-2016b, prst-2017a, prst-2017b, arachchige2019robust}, we are confident of similarly good coverage for $\rcv_Q$.  We also propose an asymptotic interval for $\rcv_M$ as well as bootstrap intervals. 

We estimate the $p$\,th quantile $x_p=G(p)=F^{-1}(p)$ by the \cite{hynd-1996}
quantile estimator $\widehat x_p=\widehat G(p)$, which is a linear combination of two adjacent order statistics.  It is readily available as the Type 8 quantile estimator on the R software \citep{R}. 

\subsection{Asymptotic confidence intervals}\label{sec:asympCI}
Let $z_\alpha =\Phi ^{-1}(\alpha )$ denote the $\alpha $ quantile of the standard normal distribution.
All our 100($1-\alpha$)\% confidence intervals for measures of relative spread $\Tf(F)$ will be of the form:
 \begin{equation}\label{eqn:ci}
    \Tf(F_n) \pm z_{1-\alpha/2}\;\widehat{\asd} (\Tf,F_n)/\sqrt n\,~,
 \end{equation}
 where $\Tf(F_n)$ is the estimator of $\Tf(F)$ and $\widehat{\asd} (\Tf,F_n)/\sqrt n\,$ is an estimate of its standard deviation (standard error) based on the sample.  The actual coverage probability of this estimator depends on how quickly the distribution of $\Tf(F_n)$ approaches normality, as well as the rate of convergence of $\Tf(F_n)$ to $\Tf(F)$ and $\widehat{\asd} (\Tf,F_n)$ to $\asd (\Tf,F).$

In constructing the interval estimators for the ratios, due to improved statistical performance such as quicker convergence to normality, it is common to first construct the interval for the log-transformed ratio followed by exponentiation to return to the original ratio scale. Let $W(F)=\ln[\Tf(F)]$ then, using the Delta Method \citep[e.g. Ch.3 of][]{Das-2008},
\begin{equation}
\asv(W,F) \doteq \frac{1}{[\Tf(F)]^2} \ \asv(\Tf,F)~.
\end{equation}
Then $\widehat{\asd} (W,F_n)= \{\widehat{\asv} (W,F_n)\}^{1/2}$, where $\widehat{\asv} (W,F_n))$ is an estimate of the asymptotic variance, enables one to construct the confidence interval for $W(F)$, which is based on the asymptotic normality of $W(F_n)$, before exponentiating to the original scale.

\subsubsection{Confidence interval for CV}\label{sec:CI_cv}
A $(1-\alpha)\times 100$\% confidence interval for the CV, which is based on the asymptotic normality of $\wcv$ when the first four moments of $F$ exist is
\begin{equation}\label{eqn:cicv}
[L,U]_{\cv } \equiv  \exp\left[\ln{(\scv)}\pm z_{1-\alpha/2}\;\frac {\widehat{\asd}{(\cvf, F_n)}}{\scv\sqrt n\,}\right]~
\end{equation}
and later we define this confidence interval method as \lq \lq Delta CV \rq\rq in our simulation study.  The ASV for the CV estimator is given in \eqref{eqn:asymvarCV} and to obtain our asymptotic standard error we replace the population CV, $\sigma$ and $\mu$ with $\scv$, sample standard deviation $s$ and sample mean $\overline{x}$ respectively.  To estimate $\mu_j$ (the $j$th central moment) we use $n^{-1}\sum^n_{i=1}(x_i-\overline{x})^j$.
\subsubsection{Confidence interval for $\rcv_Q$}\label{sec:CI_RCV_Q}
A large-sample confidence interval for $\rcv_Q=0.75\,\iqr /m$ is in terms of the estimate
$\srcv_Q=0.75(\widehat x_{0.75}-\widehat x_{0.25})/\widehat x_{0.5}$
\begin{equation}\label{eqn:ciRCV_Q}
    [L,U]_{\rcv_Q} = \exp \left [\ln{(\srcv_Q)}\pm z_{1-\alpha /2}\, \frac {\widehat{\asd} (\rcvf_Q, F_n)}{\srcv_Q\sqrt n\,}\right]~.
\end{equation}
The $\asv(\rcvf_Q, F))$ is given in Theorem \ref{th:ASV_RCVQ} and to obtain $\widehat{\asd} (\rcvf_Q, F_n)=\sqrt{\widehat{\asv} (\rcvf_Q,F_n)}$,  one needs to replace each $x_p$ by $\widehat x_p$  and each $g(p)$ by $\widehat{g}(p)$.  For $\widehat{g}(p)$, we use a kernel density estimator with the \cite{epan-1969} kernel and optimal bandwidth using the quantile optimality ratio of \cite{prst-2016a}.  

\subsubsection{Confidence interval for $\rcv_M$}\label{sec:CI_RCV_M}

A large-sample confidence interval for $\rcv_M=1.4826\,\mad /m$  is in terms of $\srcv_M=1.4826\,\widehat{\mad} /\widehat x_{0.5}$,
\begin{equation}\label{eqn:ciRCV_M}
    [L,U]_{\rcv_M} = \exp \left [\ln{(\srcv_M )}\pm z_{1-\alpha /2}\, \frac {\widehat{\asd} (\rcvf_M,F_n)}{\srcv_M\sqrt n\,}\right ]~ .
\end{equation}
Estimation of the MAD is trivial, requiring only routine coding if functionality is not already available (i.e. it is simply the median of the ordered absolute differences of the $x_i$s from the sample median).  We also need to estimate $\rho_1$, $\rho_2$ and $\rho_{12}$ in \eqref{eq:ASV_RCVM} and a simple approach using readily available software is use the FKML parameterization \citep{fmkl-1988} of the Generalized Lambda Distribution (GLD).  Defined in terms of its quantile function
$$Q(p)=\lambda_1+\frac{1}{\lambda_2}\left(\frac{p^{\lambda_3} - 1}{\lambda_3}-\frac{(1-p)^{\lambda_4} - 1}{\lambda_4}\right),$$
where $\lambda_i$ $(i=1,\ldots,4)$ are location, inverse scale and two shape parameters, the GLD can approximate a very wide range of probability distributions \citep[e.g.][]{kardud-2000, dedduwakumara2019simple}.  To do so we use the method of moments estimators and density and quantile functions for the GLD in R \texttt{gld} package \citep{gld-king}.  It is then simple to estimate  $\rho_1$, $\rho_2$ and $\rho_{12}$ using the quantile and density functions with the estimated GLD parameters and the estimated MAD.

Additional to the asymptotic interval above, we also consider two bootstrap confidence intervals. 

\subsubsection*{Non-parametric bootstrap} \label{sec:Non_para boot}
A non-parametric bootstrap re-samples $n$ observations with replacement from the sample and estimates the MAD.  This is repeated $B$ times and let $\widehat{\mad}^i$ $(i=1,\ldots,B)$ denote the $i$th estimated MAD.  The lower and upper bounds for the 95\% bootstrap interval is then the 0.025 and 0.975 quantiles of the estimated $\widehat{\mad}^i$s.

\subsubsection*{Parametric bootstrap} \label{sec:Para boot}
The parametric bootstrap interval is obtained in the same way as the non-parametric bootstrap with the exception that the sampling is done from a nominated, or estimated, density function.  In this case, we use the estimated density from the FKML GLD as described above for the asymptotic interval. This is called the Generalized Bootstrap by  \cite{dudewicz1992generalized} who also uses the GLD, albeit with a different parameterization, as one example.

\subsection{Confidence intervals for comparing two relative spreads}\label{sec:CI_diff}
When data from two independent groups are available, it is straightforward to obtain interval estimators for the comparison of relative spread for each group.  Given that empirical evidence suggests excellent coverage can be achieved in the single sample case by using a log transformation, we propose to use the log ratio of two independent relative spread estimators with a back exponentiation to the ratio scale.  For example, an interval estimator for $\rcv_{M,1}/\rcv_{M,2}$ where $\rcv_{M,1}$ and $\rcv_{M,2}$ are the relative MAD-based spread for independent populations, is, where for simplicity $\widehat{r}=\srcv_{M,1}/\srcv_{M,2}$,
\begin{equation}
\exp \left[\ln(\widehat{r})\pm z_{1-\alpha /2}\,\left\{ \frac {\widehat{\asd} \left(\rcvf_{M,1},F_n\right)}{\srcv_{M,1}\sqrt n_1\,}+\frac {\widehat{\asd} \left(\rcvf_{M,2}, F_n\right)}{\srcv_{M,2}\sqrt n_2\,}\right\}\right]
   ~,
   \label{eqn:ciRCV_M_diff}
\end{equation}
where $n_1$ and $n_2$ are the sample sizes for simple random samples from the populations and where the estimates and asymptotic standard errors can be found as above for the single sample setting.

\section{Simulations and Examples} \label{sec:simuEx}
\subsection{Simulations}\label{sec:simu}
Firstly, a simulation study was conducted to compare the performance of the interval estimator of $\rcv_Q$ and asymptotic CV interval given in \ref{sec:asympCI} with the methods given in Section \ref{sec:selected_methods} using coverage probability and width as performance measures. We have selected normal (N), log normal (LN), exponential (EXP), chi-square ($\chi^2$) and Pareto (PAR) distributions with different parameter choices and with sample sizes $n =\{50, 100, 200, 500, 1000\}$. 10,000 simulation trials were used.

\begin{landscape}
\begin{table}[htbp]
  \centering
  \scriptsize
  \caption{Simulated coverage probabilities (and widths) for 95\% confidence interval estimators for $\rcv_Q$, Delta $\cv$ and the intervals for $\cv$ described in Section \ref{sec:selected_methods}. (* median widths reported due to excessively large average widths after back-exponentiation.)} 
   \resizebox{\linewidth}{!}{%
    \begin{tabular}{rlccccccc}
    \toprule
    \multicolumn{1}{c}{Sample } & \multicolumn{1}{c}{Distribution} & Panich & Med   & Med   & Gulhar & Inverse & Delta & $\rcv_Q$ \\
    \multicolumn{1}{c}{Size(n)} &       &       &  Mill & MMcK  & Method & Method  &  CV   &    \\
    \midrule
    50    & N(5, 1) & 0.927(0.08) & 0.937(0.08) & 0.941(0.08) & 0.943(0.08) & 0.838(0.06) & 0.929(0.08) & 0.979(0.16) \\
          & LN(0, 1) & 0.688(0.97) & 0.817(1.03) & 0.803(1.07) & 0.508(0.48) & 0.808(4.85) & 0.997(7.81*) & 0.983(1.30) \\
          & EXP(1) & 0.965(0.78) & 0.978(0.73) & 0.981(0.88) & 0.887(0.40) & 0.992(3.54) & 0.997(0.68) & 0.985(1.30) \\
          & Chi(5) & 0.954(0.34) & 0.971(0.35) & 0.966(0.36) & 0.918(0.26) & 0.999(0.76) & 0.959(0.33) & 0.977(0.58) \\
          & PAR(1, 4) & 0.746(1.12) & 0.866(1.19) & 0.836(1.22) & 0.552(0.52) & 0.720(2.97) & 1.000(3.57E+9*) & 0.985(1.70) \\
    \midrule
    100   & N(5, 1) & 0.938(0.06) & 0.949(0.06) & 0.948(0.06) & 0.943(0.06) & 0.900(0.05) & 0.938(0.06) & 0.978(0.11) \\
          & LN(0, 1) & 0.755(0.85) & 0.842(0.77) & 0.867(0.96) & 0.453(0.35) & 0.926(2.69) & 0.980(5.64) & 0.975(0.82) \\
          & EXP(1) & 0.979(0.55) & 0.988(0.52) & 0.991(0.62) & 0.863(0.28) & 1.000(2.06) & 0.983(0.43) & 0.971(0.84) \\
          & Chi(5) & 0.966(0.24) & 0.961(0.34) & 0.975(0.26) & 0.909(0.18) & 1.000(0.59) & 0.953(0.22) & 0.971(0.39) \\
          & PAR(1, 4) & 0.812(0.99) & 0.887(0.88) & 0.914(1.11) & 0.471(0.37) & 0.890(1.79) & 1.000(1.19E+6*) & 0.978(1.07) \\
    \midrule
    200   & N(5, 1) & 0.947(0.04) & 0.946(0.04) & 0.945(0.04) & 0.940(0.04) & 0.955(0.04) & 0.942(0.04) & 0.979(0.08) \\
          & LN(0, 1) & 0.783(0.67) & 0.828(0.56) & 0.892(0.75) & 0.404(0.25) & 0.979(2.46) & 0.970(2.30) & 0.967(0.55) \\
          & EXP(1) & 0.987(0.39) & 0.988(0.37) & 0.997(0.43) & 0.850(0.20) & 1.000(1.44) & 0.974(0.29) & 0.966(0.57) \\
          & Chi(5) & 0.976(0.17) & 0.967(0.17) & 0.978(0.18) & 0.911(0.13) & 1.000(0.47) & 0.955(0.15) & 0.968(0.27) \\
          & PAR(1, 4) & 0.822(0.78) & 0.871(0.65) & 0.929(0.87) & 0.422(0.27) & 0.970(4.20) & 0.999(1.16E+4*) & 0.969(0.71) \\
    \midrule
    500   & N(5, 1) & 0.944(0.03) & 0.949(0.03) & 0.950(0.03) & 0.944(0.02) & 0.987(0.03) & 0.950(0.03) & 0.967(0.05) \\
          & LN(0, 1) & 0.792(0.44) & 0.782(0.36) & 0.923(0.49) & 0.360(0.16) & 0.998(2.21) & 0.965(1.21) & 0.961(0.33) \\
          & EXP(1) & 0.991(0.25) & 0.960(0.23) & 0.994(0.27) & 0.841(0.12) & 1.000(1.00) & 0.959(0.18) & 0.960(0.35) \\
          & Chi(5) & 0.976(0.11) & 0.956(0.11) & 0.966(0.11) & 0.914(0.08) & 1.000(0.36) & 0.951(0.09) & 0.960(0.17) \\
          & PAR(1, 4) & 0.833(0.52) & 0.828(0.42) & 0.952(0.58) & 0.368(0.17) & 0.995(1.85) & 0.999(291.02*) & 0.963(0.43) \\
    \midrule
    1000  & N(5, 1) & 0.952(0.02) & 0.949(0.02) & 0.951(0.02) & 0.943(0.02) & 0.997(0.03) & 0.954(0.02) & 0.960(0.03) \\
          & LN(0, 1) & 0.751(0.31) & 0.739(0.26) & 0.874(0.35) & 0.336(0.11) & 0.999(1.51) & 0.959(0.81) & 0.959(0.23) \\
          & EXP(1) & 0.992(0.18) & 0.884(0.16) & 0.964(0.19) & 0.834(0.09) & 1.000(0.79) & 0.955(0.13) & 0.958(0.24) \\
          & Chi(5) & 0.979(0.08) & 0.928(0.08) & 0.950(0.08) & 0.906(0.06) & 1.000(0.29) & 0.949(0.07) & 0.958(0.12) \\
          & PAR(1, 4) & 0.797(0.37) & 0.794(0.30) & 0.923(0.41) & 0.339(0.12) & 0.998(1.65) & 0.998(52.74*) & 0.956(0.30) \\
    \bottomrule
    \end{tabular}}%
  \label{table:Simulatios_RCV}%
\end{table}%
\end{landscape}

\begin{table}[htbp]
  \centering
  \scriptsize
  \caption{Simulated Coverage probabilities (and widths) for 95$\%$ bootstrap (non-parametric and parametric) confidence interval estimators for $\rcv_M$}
  \resizebox{\linewidth}{!}{%
    \begin{tabular}{cllll}
    \toprule
    \multicolumn{1}{l}{Sample} & \multicolumn{1}{c}{\multirow{2}[2]{*}{Distribution}} & \multicolumn{2}{c}{Method} &  \\
    \multicolumn{1}{l}{size(n)} &       & Non-parametric & Parametric & Asymptotic \\
   \midrule
    \multirow{5}[2]{*}{50} & N(5, 1) & 0.9740(0.141)  & 0.9616(0.131) & 0.9525(0.134) \\
          & LN(0, 1) & 0.9772(0.479)  &  0.9839(0.441) & 0.9665(0.524) \\
          & EXP(1) & 0.9758(0.565) &  0.9893(0.508) & 0.9719(0.601) \\
          & Chi(5)  &  0.9763(0.421)  &  0.9840(0.394) & 0.9557(0.413) \\
          &  PAR(1, 4) &  0.9777(0.549)  &  0.9874(0.493) & 0.9751(0.619) \\
    \midrule
    \multirow{5}[2]{*}{100} & N(5, 1) &  0.9759(0.099)  &  0.9795(0.093)  & 0.9493(0.094) \\
          & LN(0, 1) &  0.9749(0.337)  &  0.9859(0.327)  & 0.9673(0.370) \\
          & EXP(1) &  0.9762(0.402) & 0.9946(0.374) & 0.9648(0.411) \\
          & Chi(5)  &  0.9738(0.296)  &  0.9776(0.284) & 0.9588(0.291) \\
          &  PAR(1, 4) &  0.9748(0.389)  &  0.9933(0.362)  & 0.9697(0.414) \\
    \midrule
    \multirow{5}[2]{*}{200} & N(5, 1) & 0.9725(0.069)  &  0.9826(0.066) & 0.9520(0.066) \\
          & LN(0, 1) & 0.9724(0.235) &  0.9688(0.236)  & 0.9726(0.265) \\
          & EXP(1) &  0.9720(0.282)  &  0.9965(0.270) & 0.9591(0.287) \\
          & Chi(5)  &  0.9704(0.207)  &  0.9848(0.201) & 0.9576(0.205) \\
          &  PAR(1, 4) &  0.9729(0.272) &  0.9903(0.261)  & 0.9681(0.283) \\
    \midrule
    \multirow{5}[2]{*}{500} & N(5, 1) &  0.9644(0.043) &  0.9851(0.042)  & 0.9505(0.042) \\
          & LN(0, 1) & 0.9668(0.147) &  0.9257(0.150) & 0.9757(0.169) \\
          & EXP(1) &  0.9624(0.177)  &  0.9962(0.173)  & 0.9564(0.180) \\
          & Chi(5)  & 0.9678(0.129) &  0.9877(0.127)  & 0.9574(0.129) \\
          &  PAR(1, 4) &  0.9681(0.171)  & 0.9570(0.167)  & 0.9635(0.176) \\
    \midrule
    \multirow{5}[2]{*}{1000} & N(5, 1) &  0.9582(0.030) &  0.9861(0.029) & 0.9495(0.030) \\
          & LN(0, 1) & 0.9616(0.103)  & 0.8247(0.106)  & 0.9793(0.120) \\
          & EXP(1) &  0.9612(0.124)  &  0.9757(0.123) & 0.9569(0.128) \\
          & Chi(5)  &  0.9640(0.091)  &  0.9834(0.090) & 0.9571(0.092) \\
          &  PAR(1, 4) &  0.9606(0.119)  &  0.8029(0.118) & 0.9621(0.124) \\
    \bottomrule
    \end{tabular}}%
  \label{table:Simulatios_RCV_M}%
\end{table}%

In Table \ref{table:Simulatios_RCV} we provide the simulation results for the CV and RCV$_Q$ intervals.  For simplicity, the RCV$_M$ results follow in Table \ref{table:Simulatios_RCV_M} where the bootstrap and asymptotic intervals are compared. From Table \ref{table:Simulatios_RCV}, the Panich, Med Mill and Gulhar interval estimators for the CV perform really well for the normal distribution and when the sample size increases coverage reach to the nominal coverage. However, coverages was typically below nominal for skewed distributions pointing to unreliable performance of the estimators.  The Delta CV interval of \eqref{sec:CI_cv} provides improved coverage and close to nominal when the sample size increases, with the exception for the PAR(5,1) distribution for which the CV is undefined. The interval estimator for $\rcv_Q$ was conservative being slightly above nominal for these simulations.  The asymptotic interval for $\rcv_M$ (Table \ref{table:Simulatios_RCV_M}) provide excellent coverage, even for $n=50$ and all distributions considered.  With notable narrower intervals and very good coverage, the use of $\rcv_M$ and associated asymptotic interval estimators using estimated GLD functions are practically enticing.  However, there does not appear to be a benefit for using a bootstrap approach where coverage was typically more conservative.

\subsubsection{A Shiny web application for the performance comparisons of the intervals} \label{sect:Shiny}

For further comparisons, we have developed a Shiny \citep{shiny} web application that readers can use to run the simulations with different parameter choices.  This can be found at \url{https://lukeprendergast.shinyapps.io/Robust_CV/}.  The user can change the distribution, parameters, sample size, probability and the number of trials according to their choices.  Once the desired options are selected, the  \lq Run Simulation \rq button can be pressed and the relevant estimates, coverage probability (cp) and the average width of the confidence interval (w) will be calculated according to their input choices. In addition to that in the bottom right hand corner of the web page it will shows the time taken to run the each simulation.

\subsection{Examples}\label{sec:Ex}

We have selected two different data sets, which are named as doctor visits data and Melbourne house price data  to apply our findings to real world data.

\subsubsection{Doctor visits data}\label{sec:Dovisits}

We selected the doctor visits data set used in \cite{heritier2009robust} to apply our findings to a real world problem. The doctor visits data is a subsample of 3066 individuals of the AHEAD cohort (born before 1924) for wave 6 (year 2002) from the Health and Retirement Study (HRS) which surveys more than 22,000 Americans over the age of 50 every 2 years. We grouped this data in to two groups by taking the gender as the grouping variable. The response variable that we were interested is the number of doctor visits. Table \ref{tab:dovisits} provides summary statistics of the response variable for the two gender groups.

\begin{table}[h!t]
  \centering
  \caption{Summary Statistics of number of doctor visits between Male and Female}
    \begin{tabular}{cccc}
    \toprule
    Summary & Male  & Female & Female \\
    Statistic &       &       & (without outlier) \\
    \midrule
    Sample Size & 987   & 2079  & 2078 \\
    Minimum & 0     & 0     & 0 \\
    1st Quartile & 4     & 4     & 4 \\
    Median     & 8     & 8     & 8 \\
    Mean  & 12.08 & 12.8  & 12.45 \\
    3rd Quartile   & 14    & 15    & 15 \\
    Maximum & 300   & 750   & 365 \\
    \bottomrule
    \end{tabular}%
  \label{tab:dovisits}%
\end{table}%

From Table \ref{tab:dovisits}, the summary statistics suggest that the doctor visits distributions are positively skewed which is common for count variables.  There is also a large outlier in the female group with a number of doctor visits equal to 750. We removed the outlier form the data set and again calculated the descriptive statistics for female group as shown in the $3^{rd}$ column of the above Table \ref{tab:dovisits}. The mean for the female group reduces after the removal of the outlier and the summary statistics still suggest positive skew. 

Our objective was to compare the relative spread of the number of doctor visits between males and females. We used $\cv$, $\rcv_Q$ and $\rcv_M$ to compare the relative spread of the number of doctor visits between males and females with and without an outlier.  

\begin{table}[htbp]
  \centering
  \scriptsize
  \caption{95 \% confidence interval lower bounds (LB) and upper bounds (UB) for the number of doctor visits.}
  \vspace{0.1cm}
 \resizebox{\linewidth}{!}{%
\begin{tabular}{lccc}
\toprule
    Sample & CV & RCV$_Q$ & RCV$_M$  \\
    \midrule
    Male & $(1.283,\ 2.016)$ & $(0.837,\    1.050)$ & $(0.681,\   0.807)$ \\ 
    Female & $(1.298,\    2.801)$ & $(0.943,\    1.128)$ & $(0.700,\   0.786)$ \\ 
    Female, outlier excluded & $(1.237,\ 1.746)$ & $(0.943,\    1.128)$ & $(0.699,\   0.786)$ \\ 
    \bottomrule
\end{tabular}}
  \label{tab:CI_dovisits}%
\end{table}%

Table \ref{tab:CI_dovisits} provides the confidence interval bounds of the 95 percent confidence intervals for the three measures.  The confidence interval for CV is greatly influenced by whether or not the outlier in the female data is included.  This is not the case for the interval for quantile-based measures.  Additionally, in comparison, the interval $\cv$  is wide compared to the intervals for $\rcv_Q$ and $\rcv_M$.

\subsubsection{Melbourne house price data} \label{sect:HousePrice} 
The median is the most popular summary measure used to describe housing markets. Motivated by this, we applied our measures to Melbourne house clearance data from January 2016 which  is available at \url{https://www.kaggle.com/anthonypino/melbourne-housing-market}. This data set contains suburb-wise prices for three types of houses (house, unit, townhouse). There is data for 369 suburbs and we removed the suburbs, which contain less than 10 houses sold leaving 301 suburbs.

We selected three pairs of suburbs which were considered by \citep{arachchige2019robust} to calculate the interval estimators for ratios $\cv$, $\rcv_Q$ and $\rcv_M$ to assess differences in relative spread of house prices. 

\begin{figure}[!htb]
 \centering
  \includegraphics[width=\linewidth]{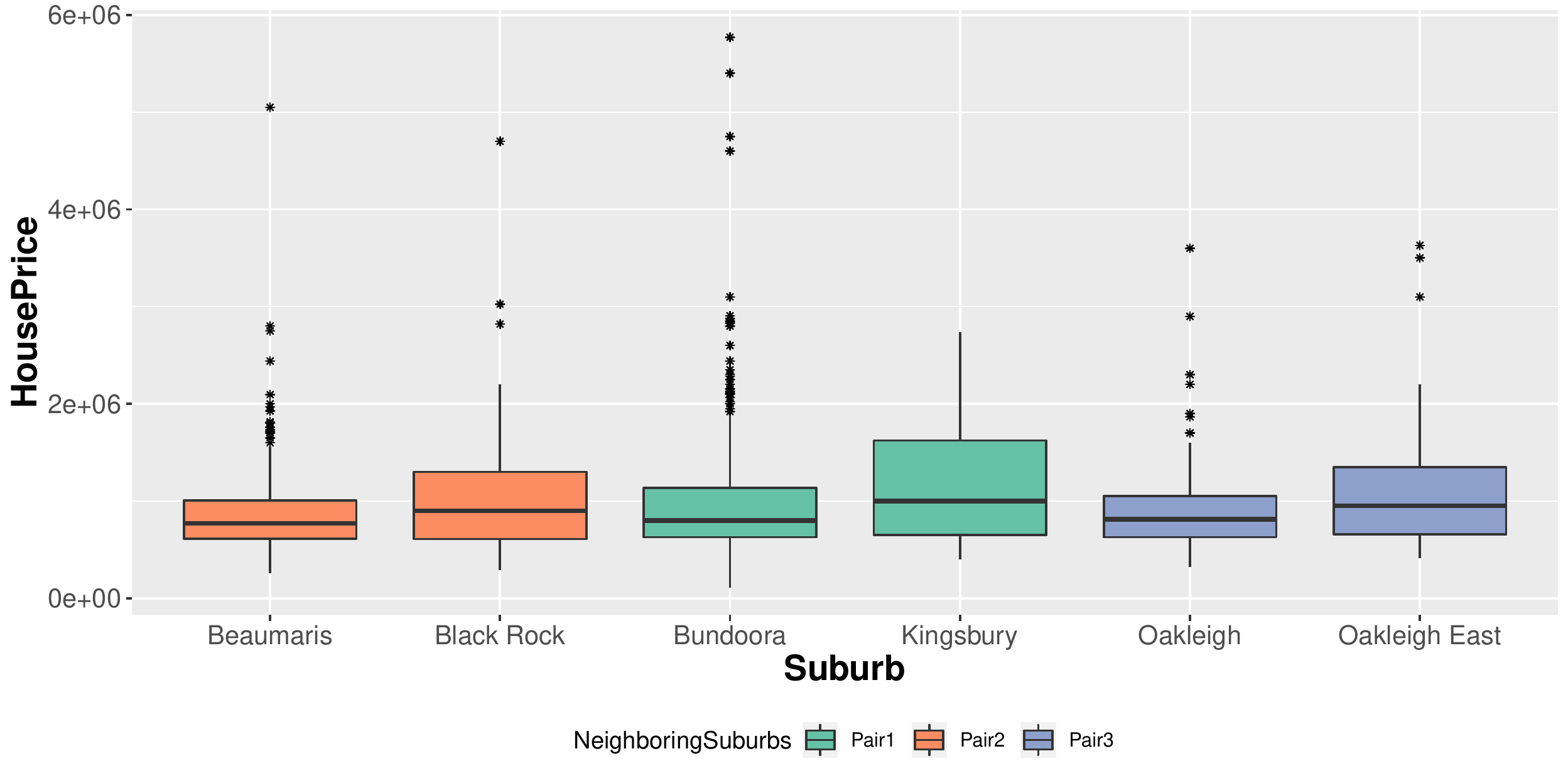}
  \caption{House price comparisons of selected three pair of neighboring suburbs}\label{fig:HP_neigboxplots}
\end{figure}

Figure \ref{fig:HP_neigboxplots} depicts there are outliers for all suburbs except for Kingsbury. Additionally, there are differences in spread for the house price distributions between each neighboring suburb.

\begin{table}[htbp]
  \centering
  \scriptsize
  \caption{95 \% confidence interval lower bounds (LB) and upper bounds (UB) for ratios of $\cv$, $\rcv_Q$ and $\rcv_M$ between neighboring suburbs house prices.}
  \resizebox{\linewidth}{!}{%
    \begin{tabular}{lcccccc}
    \toprule
    \multicolumn{1}{c}{Confidence} & \multicolumn{2}{c}{$x=$Bundoora} & \multicolumn{2}{c}{$x=$Black Rock} & \multicolumn{2}{c}{$x=$Oakleigh} \\
    \multicolumn{1}{c}{ Interval} & \multicolumn{2}{c}{$y=$Kingsbury} & \multicolumn{2}{c}{$y=$Beaumaris} & \multicolumn{2}{c}{$y=$Oakleigh East} \\
    \multicolumn{1}{c}{Method} & LB    & UB    & LB    & UB    & LB    & UB \\
    \midrule
    $\cv_x/\cv_y$    & \multicolumn{1}{r}{1.0156} & 1.6079 & \multicolumn{1}{r}{0.6525} & 1.3225 & \multicolumn{1}{r}{0.7219} & 1.3519 \\
    $\rcv_{Q_x}/\rcv_{Q_y}$    & \multicolumn{1}{r}{0.4336} & \multicolumn{1}{r}{0.9736} & \multicolumn{1}{r}{0.4844} & 0.9243 & \multicolumn{1}{r}{0.4607} & 1.0914 \\
     $\rcv_{M_x}/\rcv_{M_y}$    & \multicolumn{1}{r}{0.5392} & \multicolumn{1}{r}{1.0808} & \multicolumn{1}{r}{0.5751} & 0.9366 & \multicolumn{1}{r}{0.5286} & 1.0218 \\
    \bottomrule
    \end{tabular}}%
  \label{tab:CIdif_HPrice}%
\end{table}%
\FloatBarrier

Ratios of the measure are reported in Table \ref{tab:CIdif_HPrice} to see whether there is a difference in relative spread between suburbs.  Comparing Bundoora and Kingsbury, the measures provide different insights.  While the box plot suggests greater spread in Kingsbury, the ratio of CVs suggests otherwise having been highly influence by outliers in Bundoora.  The ratios of RCV$_Q$ and RCV$_M$ suggest greater relative spread in Kingsbury which is in better agreement with what is shown in the box plots.  For Beaumaris and Black Rock, a significant difference is not found for the CVs and the interval is wide.  However, the other intervals suggest a significant difference.  All three measures suggest there is not a significant difference in relative spread of house price between Oakleigh and Oakleigh East, although the intervals do tend to suggest that there is for RCV$_Q$ and RCV$_M$.  Overall, the intervals are narrower for the quantile-based measures having not been so greatly influence by outliers.
\section{Summary and discussion}\label{sec:summary}
We have proposed interval estimators for alternative robust measures of relative spread to the coefficient of variation.  RCV$_Q$, a scalar multiple of the interquartile range divided by median, is simple and the associated confidence intervals have very good coverage over a diverse range of distribution types.  Similarly, RCV$_M$ where the MAD is used instead of the interquartile range, interval also have excellent coverage and typically has smaller variability than the estimator for RCV$_Q$ making it a preferred candidate to be used instead of the CV.  While we also considered bootstrap interval estimators for RCV$_M$, the asymptotic Wald-type interval based on the approximate variances, and covariance between, the MAD and median achieved excellent coverage even for sample sizes as small as 50.  These robust intervals compare very favorably to the CV where coverage is typically poor when the data is not sampled from a normal distribution.  Our examples highlighted that they can provide very different insights into relative spread when compared to the CV, and the use of quantile-based measures is more easily justified when data is skewed due to difficulty interpreting the mean and variance.

\appendix

\section{Proof of Theorem \ref{th:ASV_RCVQ}}\label{app:proof_of_ASV_RCVQ}
Recall $\IF (x;\,\rho _{p,q} ,F)$ and $\IF (x;\,\rcvf_Q ,F)$ in \eqref{eq:IF_rho} and \eqref{eq:IF_RCVq} respectively. For simplicity let $\IF (x;\,\rho _{p,q} ,F)=\IF_{\rho _{p,q}}$,
$\IF (x;\,\rcvf_Q ,F)=\IF_{\rcvf_Q}$, $\IF [\Gf(\,\cdot ,p)]=\IF_{\Gf ,p}$ and $\asv\left(\Gf,F; p\right)=\asv_{\Gf, p}$. Then
\begin{equation}
\e(\IF_{\rcvf_Q}^2) = 0.75^2\,\left[ \e\left(\IF_{\rho _{3/4,1/2}}^2 \right) + \e\left(\IF_{\rho _{1/4,1/2}}^2 \right)  -  2  \e \left(\IF_{\rho _{3/4,1/2}} \IF_{\rho_{1/4,1/2}}\right)\right]. 
   \label{eq:ASV_RCVQ1}
\end{equation}
 It can be shown, 
 \begin{align}
  \e\left(\IF_{\rho_{3/4,1/2}}^2 \right)  &= \rho_{3/4,1/2}^2\ \e\left[\left( \frac {\IF _{\Gf,3/4}}{x_{3/4}}-\frac {\IF_{\Gf,1/2}}{x_{1/2}}\right)\right ]^2 \nonumber\\
   &= \frac{x_{3/4}^2}{x_{1/2}^2}\left[\frac {\e\left(\IF _{\Gf,3/4}^2\right)}{x_{3/4}^2}
          +\frac {\e\left(\IF_{\Gf,1/2}^2\right)}{x_{1/2}^2}  -  \frac{2\e\left(\IF_{\Gf,3/4} \IF_{\Gf,1/2}\right)}{x_{3/4}x_{1/2}} \right] \nonumber\\
       &= \frac{1}{x_{1/2}^2} \left[\asv_{\Gf,3/4}
         +\frac {x_{3/4}^2\asv_{\Gf,1/2}}{x_{1/2}^2}-  \frac{2x_{3/4}\e\left(\IF_{\Gf,3/4}\ \IF_{\Gf,1/2}\right)}{x_{1/2}} \right]~.
     \label{eq:rho1}
  \end{align}
 Similarly,
  \begin{align}
  \e\left(\IF_{\rho_{1/4,1/2}}^2 \right) =& \frac{1}{x_{1/2}^2} \left[\asv_{\Gf,1/4}
         +\frac {x_{1/4}^2\asv_{\Gf,1/2}}{x_{1/2}^2}  -  \frac{2x_{1/4}\e\left(\IF_{\Gf,1/4}\ \IF_{\Gf,1/2}\right)}{x_{1/2}} \right]  
 \label{eq:rho2}
  \end{align}
  and
\begin{align}\label{eq:rho12}
   \e \left(\IF_{\rho _{3/4,1/2}} \, \IF_{\rho _{1/4,1/2}}\right)
   =& \rho_{3/4,1/2}\times\rho_{1/4,1/2}\, \e \Bigg[\left( \frac {\IF_{\Gf,3/4}}{x_{3/4}}
          -\frac {\IF_{\Gf,1/2}}{x_{1/2}}\right) \nonumber\\
         & \qquad \qquad \times \left(\frac {\IF_{\Gf,1/4}}{x_{1/4}}
          -\frac {\IF_{\Gf,1/2}}{x_{1/2}}\right) \Bigg] \nonumber\\
 =& \frac{1}{x_{1/2}^2} \Bigg[ \e \left(\IF_{\Gf,3/4}\,\IF_{\Gf,1/4}\right)  -  \frac {x_{1/4}\e \left(\IF_{\Gf,3/4}\,\IF_{\Gf,1/2} \right)}{x_{1/2}}  \nonumber\\
& -  \frac {x_{3/4}\e \left(\IF_{\Gf,1/4}\,\IF_{\Gf,1/2} \right)}{x_{1/2}}  +\frac {x_{3/4}x_{1/4}\e\left(\IF_{\Gf,1/2}^2\right)}{x_{1/2}^2}\Bigg] ~. 
 \end{align}
 Substituting the above \eqref{eq:rho1}, \eqref{eq:rho2},\eqref{eq:rho12} in \eqref{eq:ASV_RCVQ1} and using $\asv\left(\Gf,F; p\right)=p(1-p)g^2(p)$ gives

  \begin{align}
   \e[\IF_{\rcvf_Q}^2]
=& \frac{0.75^2(x_{3/4}-x_{1/4})^2}{x_{1/2}^2}\,  \Bigg\{\frac{\asv_{\Gf,3/4}+\asv_{\Gf,1/4} 
  -2  \e \left(\IF_{\Gf,3/4}\,\IF_{\Gf,1/4}\right)}{(x_{3/4}-x_{1/4})^2} \nonumber\\
     &  +  \frac{\asv_{\Gf,1/2}}{x_{1/2}^2} 
     -\frac{2\, \left[\e \left(\IF_{\Gf,3/4}\,\IF_{\Gf,1/2} \right) 
   -  \e\left(\IF_{\Gf,1/4}\,\IF_{\Gf,1/2} \right)\right]}{x_{1/2}(x_{3/4}-x_{1/4})} 
  \Bigg\} \nonumber\\
  =&\frac{\rcvf^2_Q}{4}\ \Bigg\{\frac {3\left[g^2(3/4)+g^2(1/4)\right]-2\,g(1/4)g(3/4)}{4\times IQR^2}\nonumber\\
  &\qquad \qquad + \frac {g^2(1/2)}{m^2} - \frac {g(1/2)\left[g(3/4)-g(1/4)\right]} {m\times IQR} \Bigg\}~.
 \end{align}

\section{Computing the true MAD}\label{Rcode}
Computing the true value of $\mad$ is not a trivial task.  We provide an R function below that can be uses to compute true value of the $\mad$ for a user-specified distribution.
\begin{verbatim}
mad <- function(dist, param){
  # Computes the true value of the MAD for a specific  
  # distribution with desired parameter choices.
  #
  # Args:
  #   dist: The distribution whose MAD 
  #         is to be calculated.
  #   param: The parameter choices of the selected   
  #          distribution whose MAD is to be calculated. 
  #
  # Returns:
  # The true value of the MAD for a specific 
  # distribution with desired parameter choices.
  
  qf <- paste0("q", dist)
   m <- do.call(qf, c(p = 0.5, param)) # find median
  abs.x.m <- function(x, dist, param, m){
    df <- paste0("d", dist)
    do.call(df, c(x = x + m, param)) 
            + do.call(df, c(x = - x + m, param))
  }
  abs.x.m.vec <- Vectorize(abs.x.m, "x")
  
  f <- function(x, dist, param, m){
    integrate(abs.x.m.vec, lower = 0, upper = x,
        dist = dist, param = param, m = m)$value - 0.5
  }
  upper <- abs(do.call(qf, c(p = 0.75, param)) + m)
  uniroot(f, interval = c(0, upper), dist = dist,
          param = param, m = m)$root
}
mad("lnorm", list(meanlog=0, sdlog=1))
mad("exp", list(rate=1))
\end{verbatim}

\bibliographystyle{authordate4}
\bibliography{ref}

\begin{thebibliography}{}

\bibitem[\protect\citename{Andersen, }2008]{andersen2008modern}
{\sc Andersen, R.} 2008.
\newblock {\em Modern methods for robust regression}.
\newblock Sage.

\bibitem[\protect\citename{Arachchige {\em et~al.\ }\relax,
  }2019]{arachchige2019robust}
{\sc Arachchige, C.~NPG, Cairns, M., \& Prendergast, L.~A.} 2019.
\newblock Interval estimators for ratios of independent quantiles and
  interquantile ranges.
\newblock {\em {Commun. Stat. B-Simul. (accepted, June)}}.

\bibitem[\protect\citename{Atkinson, }1970]{atkinson1970measurement}
{\sc Atkinson, A.~B.} 1970.
\newblock On the measurement of inequality.
\newblock {\em {J. Econ. Theor.}}, {\bf 2}(3), 244--263.

\bibitem[\protect\citename{Bonett, }2006]{bonett2006confidence}
{\sc Bonett, D.~G.} 2006.
\newblock Confidence interval for a coefficient of quartile variation.
\newblock {\em {Comput. Stat. Data An.}}, {\bf 50}(11), 2953--2957.

\bibitem[\protect\citename{Bonett \& Seier, }2005]{bonett2005confidence}
{\sc Bonett, D.~G., \& Seier, E.} 2005.
\newblock Confidence interval for a coefficient of dispersion in nonnormal
  distributions.
\newblock {\em {Biometrical J.}}, {\bf 47}(1), 144--148.

\bibitem[\protect\citename{Bulent \& Hamza, }2018]{bulent2018bootstrap}
{\sc Bulent, A., \& Hamza, G.} 2018.
\newblock Bootstrap confidence intervals for the coefficient of quartile
  variation.
\newblock {\em {Commun. Stat. B-Simul.}}, {\bf In Press}, 1--9.

\bibitem[\protect\citename{Chang {\em et~al.\ }\relax, }2017]{shiny}
{\sc Chang, W., Cheng, J., Allaire, J.~J., Xie, Y, \& McPherson, J.} 2017.
\newblock {\em shiny: Web application framework for r}.
\newblock R package version 1.0.5.

\bibitem[\protect\citename{Chen \& Fleisher, }1996]{chen1996regional}
{\sc Chen, J., \& Fleisher, B.~M.} 1996.
\newblock Regional income inequality and economic growth in china.
\newblock {\em {J. Comp. Econ.}}, {\bf 22}(2), 141--164.

\bibitem[\protect\citename{Cole {\em et~al.\ }\relax,
  }2000]{cole2000establishing}
{\sc Cole, T.~J., Bellizzi, M.~C., Flegal, K.~M., \& Dietz, W.~H.} 2000.
\newblock Establishing a standard definition for child overweight and obesity
  worldwide: international survey.
\newblock {\em {BMJ Brit. Med. J.}}, {\bf 320}(7244), 1240.

\bibitem[\protect\citename{DasGupta, }2006]{Das-2008}
{\sc DasGupta, A.} 2006.
\newblock {\em Asymptotic {T}heory of {S}tatistics and {P}robability}.
\newblock New York, NY: Springer.

\bibitem[\protect\citename{David, }1981]{david-1981}
{\sc David, H.A.} 1981.
\newblock {\em Order {S}tatistics}.
\newblock New York: John Wiley \& Sons.

\bibitem[\protect\citename{Dedduwakumara {\em et~al.\ }\relax,
  }2019]{dedduwakumara2019simple}
{\sc Dedduwakumara, D.~S., Prendergast, L.~A., \& Staudte, R.~G.} 2019.
\newblock A simple and efficient method for finding the closest generalized
  lambda distribution to a specific model.
\newblock {\em {Cogent Math. Stat.}}, {\bf 6}, 1--11.

\bibitem[\protect\citename{{D}evelopment~{C}ore {T}eam, }2008]{R}
{\sc {D}evelopment~{C}ore {T}eam, R}. 2008.
\newblock {\em R: A language and environment for statistical computing}.
\newblock R {F}oundation for {S}tatistical {C}omputing, Vienna, Austria.
\newblock {ISBN} 3-900051-07-0.

\bibitem[\protect\citename{Dudewicz, }1992]{dudewicz1992generalized}
{\sc Dudewicz, E.~J.} 1992.
\newblock The generalized bootstrap.
\newblock {\em Pages  31--37 of:} {\em Bootstrapping and related techniques}.
\newblock Springer.

\bibitem[\protect\citename{Epanechnikov, }1969]{epan-1969}
{\sc Epanechnikov, V.~A.} 1969.
\newblock Nonparametric estimation of a multivariate probability density.
\newblock {\em {Theor. Probab. Appl.}}, {\bf 14}, 153--158.

\bibitem[\protect\citename{Falk, }1997]{falk1997asymptotic}
{\sc Falk, M.} 1997.
\newblock Asymptotic independence of median and mad.
\newblock {\em {Stat. Probabil. Lett.}}, {\bf 34}(4), 341--345.

\bibitem[\protect\citename{Freimer {\em et~al.\ }\relax, }1988]{fmkl-1988}
{\sc Freimer, M., Mudholkar, G.~S., Kollia, G., \& Lin, C.~T.} 1988.
\newblock A study of the generalized {T}ukey lambda family.
\newblock {\em {Comm. Stat. A-Theor.}}, {\bf 17}, 3547--3567.

\bibitem[\protect\citename{Gastwirth, }1982]{gastwirth1982statistical}
{\sc Gastwirth, J.~L.} 1982.
\newblock Statistical properties of a measure of tax assessment uniformity.
\newblock {\em {J. Stat. Plan. Infer.}}, {\bf 6}(1), 1--12.

\bibitem[\protect\citename{Gong \& Li, }1999]{gong1999relationship}
{\sc Gong, J., \& Li, Y.} 1999.
\newblock Relationship between the estimated weibull modulus and the
  coefficient of variation of the measured strength for ceramics.
\newblock {\em {J. Am. Ceram. Soc.}}, {\bf 82}(2), 449--452.

\bibitem[\protect\citename{Groeneveld, }2011]{groen-2011}
{\sc Groeneveld, R.~A.} 2011.
\newblock Influence functions for the coefficient of variation, its inverse,
  and cv comparisons.
\newblock {\em {Comm. Stat. A-Theor.}}, {\bf 40}(23), 4139--4150.

\bibitem[\protect\citename{Gulhar {\em et~al.\ }\relax,
  }2012]{gulhar2012comparison}
{\sc Gulhar, M., Kibria, G., Albatineh, A., \& Ahmed, N.~U.} 2012.
\newblock A comparison of some confidence intervals for estimating the
  population coefficient of variation: a simulation study.
\newblock {\em Sort}, {\bf 36}(1).

\bibitem[\protect\citename{Hamer {\em et~al.\ }\relax, }1995]{hamer1995new}
{\sc Hamer, A.~J., Strachan, J.~R., Black, M.~M., Ibbotson, C., \& Elson,
  R.~A.} 1995.
\newblock A new method of comparative bone strength measurement.
\newblock {\em {J. Med. Eng. Technol.}}, {\bf 19}(1), 1--5.

\bibitem[\protect\citename{Hampel, }1974]{hamp-1974}
{\sc Hampel, F.~R.} 1974.
\newblock The influence curve and its role in robust estimation.
\newblock {\em {J. Am. Stat. Assoc.}}, {\bf 69}, 383--393.

\bibitem[\protect\citename{Hampel {\em et~al.\ }\relax, }1986]{HRRS86}
{\sc Hampel, F.~R., Ronchetti, E.~M., Rousseeuw, P.~J., \& Stahel, W.~A.} 1986.
\newblock {\em Robust {S}tatistics: The {A}pproach {B}ased on {I}nfluence
  {F}unctions}.
\newblock New {Y}ork: John {W}iley and {S}ons.

\bibitem[\protect\citename{Heritier {\em et~al.\ }\relax,
  }2009]{heritier2009robust}
{\sc Heritier, S., Cantoni, E., Copt, S., \& Victoria-Feser, M.-P.} 2009.
\newblock {\em Robust methods in biostatistics}.
\newblock  Vol. 825.
\newblock John Wiley \& Sons.

\bibitem[\protect\citename{Huber, }1981]{huber1981robust}
{\sc Huber, P.~J.} 1981.
\newblock {\em Robust statistics}.
\newblock Wiley.

\bibitem[\protect\citename{Hyndman \& Fan, }1996]{hynd-1996}
{\sc Hyndman, R.~J., \& Fan, Y.} 1996.
\newblock Sample quantiles in statistical packages.
\newblock {\em {Am. Stat.}}, {\bf 50}, 361--365.

\bibitem[\protect\citename{Karian \& Dudewicz, }2000]{kardud-2000}
{\sc Karian, Z.~A., \& Dudewicz, E.~J.} 2000.
\newblock {\em Fitting {S}tatistical {D}istributions: the generalized lambda
  distribution and generalized bootstrap methods}.
\newblock Chapman and {H}all.

\bibitem[\protect\citename{King {\em et~al.\ }\relax, }2016]{gld-king}
{\sc King, R., Dean, B., \& Klinke, S.} 2016.
\newblock {\em {gld: Estimation and Use of the Generalised (Tukey) Lambda
  Distribution}}.
\newblock R package version 2.4.1.

\bibitem[\protect\citename{Lovitt \& Holtzclaw, }1929]{lovitt1929statistics}
{\sc Lovitt, W.~V., \& Holtzclaw, H.~F.} 1929.
\newblock {\em Statistics}.
\newblock Prentice-Hall, Incorporated.

\bibitem[\protect\citename{McKay, }1932]{mckay1932distribution}
{\sc McKay, A.~T.} 1932.
\newblock Distribution of the coefficient of variation and the extended" t"
  distribution.
\newblock {\em {J. R. Stat. Soc.}}, {\bf 95}(4), 695--698.

\bibitem[\protect\citename{Miller, }1991]{edward1991asymptotic}
{\sc Miller, E.~G.} 1991.
\newblock Asymptotic test statistics for coefficients of variation.
\newblock {\em {Comm. Stat. A-Theor.}}, {\bf 20}(10), 3351--3363.

\bibitem[\protect\citename{Miller \& Karson, }1977]{miller1977testing}
{\sc Miller, E.~G., \& Karson, M.~J.} 1977.
\newblock Testing equality of two coefficients of variation.
\newblock {\em Pages  278--283 of:} {\em {Am. Stat. Assoc.: Proceedings of the
  Bus. Econ. Section, Part I}},  vol. 95.

\bibitem[\protect\citename{Panichkitkosolkul,
  }2009]{panichkitkosolkul2009improved}
{\sc Panichkitkosolkul, W.} 2009.
\newblock Improved confidence intervals for a coefficient of variation of a
  normal distribution.
\newblock {\em {Thailand Statistician}}, {\bf 7}(2), 193--199.

\bibitem[\protect\citename{Prendergast \& Staudte, }2016a]{prst-2016a}
{\sc Prendergast, L.~A., \& Staudte, R.~G.} 2016a.
\newblock {Exploiting the Quantile Optimality Ratio in finding Confidence
  Intervals for Quantiles}.
\newblock {\em Stat}, {\bf 5}, 70--81.

\bibitem[\protect\citename{Prendergast \& Staudte, }2016b]{prst-2016b}
{\sc Prendergast, L.~A., \& Staudte, R.~G.} 2016b.
\newblock Quantile versions of the {L}orenz curve.
\newblock {\em {Electron. J. Stat.}}, {\bf 10}(2), 1896--1926.

\bibitem[\protect\citename{Prendergast \& Staudte, }2017a]{prst-2017a}
{\sc Prendergast, L.A., \& Staudte, R.G.} 2017a.
\newblock When large n is not enough---distribution-free interval estimators
  for ratios of quantiles.
\newblock {\em {J. Econ. Inequal.}},  1--17.
\newblock doi:10.1007/s10888-017-9347-9.

\bibitem[\protect\citename{Prendergast \& Staudte, }2017b]{prst-2017b}
{\sc Prendergast, L.A., \& Staudte, R.G.} 2017b.
\newblock A simple and effective inequality measure.
\newblock {\em {Am. Stat.}}

\bibitem[\protect\citename{Reed {\em et~al.\ }\relax, }2002]{reed2002use}
{\sc Reed, G.~F., Lynn, F., \& Meade, B.~D.} 2002.
\newblock Use of coefficient of variation in assessing variability of
  quantitative assays.
\newblock {\em {Clin. Diagn. Lab. Immun.}}, {\bf 9}(6), 1235--1239.

\bibitem[\protect\citename{Reimann {\em et~al.\ }\relax,
  }2008]{reimann2008statistical}
{\sc Reimann, C., Filzmoser, P., Garrett, R.~G., \& Dutter, R.} 2008.
\newblock {\em Statistical data analysis explained: applied environmental
  statistics with r}.

\bibitem[\protect\citename{Shapiro, }2005]{shapiro2005practical}
{\sc Shapiro, H.~M.} 2005.
\newblock {\em Practical flow cytometry}.
\newblock John Wiley \& Sons.

\bibitem[\protect\citename{Sharma \& Krishna, }1994]{sharma1994asymptotic}
{\sc Sharma, K.~K., \& Krishna, H.} 1994.
\newblock Asymptotic sampling distribution of inverse coefficient-of-variation
  and its applications.
\newblock {\em {IEEE T. Reliab.}}, {\bf 43}(4), 630--633.

\bibitem[\protect\citename{Staudte \& Sheather, }1990]{S-S-1990}
{\sc Staudte, R.~G., \& Sheather, S.~J.} 1990.
\newblock {\em Robust {E}stimation and {T}esting}.
\newblock New York: Wiley.

\bibitem[\protect\citename{Vangel, }1996]{vangel1996confidence}
{\sc Vangel, M.~G.} 1996.
\newblock Confidence intervals for a normal coefficient of variation.
\newblock {\em {Am. Stat.}}, {\bf 50}(1), 21--26.

\bibitem[\protect\citename{Varmuza \& Filzmoser,
  }2009]{varmuza2009introduction}
{\sc Varmuza, K., \& Filzmoser, P.} 2009.
\newblock {\em Introduction to multivariate statistical analysis in
  chemometrics}.
\newblock CRC press.

\bibitem[\protect\citename{Wilcox, }2011]{wilcox2011introduction}
{\sc Wilcox, R.~R.} 2011.
\newblock {\em Introduction to robust estimation and hypothesis testing}.
\newblock Academic press.

\end{thebibliography}
\end{document}